\documentclass[draft,12pt]{amsart}
\usepackage{amssymb}

\theoremstyle{plain}
\newtheorem{theorem}{Theorem}[section]
\newtheorem{corollary}[theorem]{Corollary}
\newtheorem{lemma}[theorem]{Lemma}
\newtheorem{proposition}[theorem]{Proposition}

\theoremstyle{definition}
\newtheorem{definition}[theorem]{Definition}
\newtheorem{remark}[theorem]{Remark}
\newtheorem{example}[theorem]{Example}
\newtheorem*{question}{Question}

\theoremstyle{remark}

\newcommand{\abs}[1]{\lvert#1\rvert}
\newcommand{\norm}[1]{\lVert#1\rVert}
\newcommand{\bigabs}[1]{\bigl\lvert#1\bigr\rvert}
\newcommand{\bignorm}[1]{\bigl\lVert#1\bigr\rVert}
\newcommand{\Bigabs}[1]{\Bigl\lvert#1\Bigr\rvert}
\newcommand{\Bignorm}[1]{\Bigl\lVert#1\Bigr\rVert}
\renewcommand{\le}{\leqslant}
\renewcommand{\ge}{\geqslant}

\newcommand{\term}[1]{{\textit{\textbf{#1}}}}

\newcommand{\eqv}[1]{\overset{#1}{\approx}}

\def\N{\mathbb N}
\def\B{\mathcal B}
\def\C{\mathbb C}
\def\L{\mathcal L}
\def\K{\mathcal K}
\def\S{\mathcal S}
\def\SS{\mathcal{SS}}
\def\FSS{\mathcal{FSS}}
\def\SP{\textrm{SP}}
\def\nSPxi{\#\bigl(\SP_{1,w}(X)/\!\!\approx_\xi\bigr)}
\def\nSP{\#\bigl(\SP_{1,w}(X)/\!\!\approx\bigr)}
\def\tr{\mathrm{Tr}}
\def\wf{\mathrm{WF}}
\def\bt{\N^{<\N}}

\def\ibbt{[\N]^{<\N}}
\def\Tm{T\bigl(m,(x_n),Y\bigr)}

\DeclareMathOperator{\dist}{dist}
\DeclareMathOperator{\supp}{supp}
\DeclareMathOperator{\subs}{Subs}

\begin{document}
\baselineskip 18pt

\title{Classes of strictly singular operators and their products}

\author[G.~Androulakis]{G.~Androulakis}
\address{Department of Mathematics,
          University of South Carolina,
          Columbia, SC 29208, USA}
\email{giorgis@math.sc.edu}

\author[P.~Dodos]{P.~Dodos}
\address{National Technical University of Athens, Faculty of Applied
          Sciences, Department of Mathematics, Zografou Campus, 157
          80, Athens. Greece}
\email{pdodos@math.ntua.gr}

\author[G.~Sirotkin]{G.~Sirotkin}
\address{Department of Mathematical Sciences,
          Northern Illinois University,
          DeKalb, IL 60115. USA}
\email{sirotkin@math.niu.edu}

\author[V.~G.~Troitsky]{V.~G.~Troitsky}
\address{Department of Mathematical and Statistical Sciences,
         University of Alberta, Edmonton, AB, T6G\,2G1. Canada.
         }
\email{vtroitsky@math.ualberta.ca}

\thanks{This project originated at SUMIRFAS conference in 2005 in
  Texas A{\&}M University. The authors wish to thank the organizers of
  the SUMIRFAS conference for their hospitality}



\date{December 18, 2006}

\begin{abstract}
  V.~D. Milman proved in~\cite{Milman:70} that the product of two
  strictly singular operators on $L_p[0,1]$ ($1\le p<\infty$) or on
  $C[0,1]$ is compact.  In this note we utilize Schreier families
  $\S_\xi$ in order to define the class of $\S_\xi$-strictly
  singular operators, and then we refine the technique of Milman to
  show that certain products of operators from this class are compact,
  under the assumption that the underlying Banach space has finitely
  many equivalence classes of Schreier-spreading sequences. Finally we
  define the class of $\S_\xi$-hereditarily indecomposable Banach
  spaces and we examine the operators on them.
\end{abstract}

\maketitle

\section{Introduction}

In this paper we extend work of V.D.~Milman \cite{Milman:70} who
showed that the product of two strictly singular (bounded linear)
operators on $L_p[0,1]$ ($1\le p < \infty$) or on $C[0,1]$ is compact.
The importance of this fundamental result of V.D.~Milman lies in the
fact that compact operators are well understood, unlike strictly
singular ones.

In the present paper we use the Schreier families $\S _{\xi}$ for $1
\le \xi < \omega_1$ which were introduced by D.~Alspach and
S.A.~Argyros \cite{AA} and we define the classes of $\S_\xi$-strictly
singular operators. These classes are increasing in $\xi$ (i.e., if
$\xi < \zeta$ then every $\S_\xi$-strictly singular operator is an
$\S_\zeta$-strictly singular operator) and they exhaust the class of
strictly singular operators defined on separable Banach spaces (i.e.,
every strictly singular operator between separable Banach spaces is
$\S_\xi$-strictly singular for some $\xi$, Theorem~\ref{t1}).  We
define the notion of Schreier spreading sequence which is closely
related to the well studied notion of spreading model. In fact every
seminormalized basic sequence has a Schreier spreading subsequence.
For $1 \le \xi < \omega_1$ we define an equivalence relation
$\approx_\xi$ on the set of weakly null Schreier spreading sequences
of a Banach space. One of the main results of the present paper, which
is a refinement of the above mentioned result of V.D.~Milman, is our
Theorem~\ref{main}.  Its statement is slightly stronger than the
following simplified version: for a Banach space $X$ and an ordinal $1
\le \xi < \omega_1$, if the number of the equivalence classes of the
weakly null spreading sequences in $X$ with respect to the equivalence
relation $\approx_\xi$ is equal to $n< \infty$, then the product of
any $n+1$ many $\S_\xi$ strictly singular operators on $X$ is compact.
Applications of this result are made to Tsirelson type spaces, Read's
space~\cite{Read} and the invariant subspace problem. Finally, for $1
\le \xi < \omega_1$ we define the notion of $\S_\xi$ hereditarily
indecomposable Banach space as a refinement of the notion of
hereditarily indecomposable (HI) Banach space which was introduced by
W.T.~Gowers and B.~Maurey~\cite{Gowers:93}. If $\xi < \zeta$ then
every $\S_\xi$-HI Banach space is an $\S_\zeta$-HI space and if $X$ is
a separable HI space then it is $\S_\xi$-HI for some $1 \le \xi <
\omega_1$ (Theorem~\ref{t1}). The study of operators on complex
$\S_\xi$-HI Banach spaces and their subspaces reveals that the
$\S_\xi$-strictly singular operators play an analogous role that
strictly singular operators play on the analysis of operators on
complex HI spaces. This indicates a potential use of HI spaces towards
the solution of the invariant subspace problem
(Corollary~\ref{hyperinv}).

We recall the definition of the \term{Schreier families} $\S _\xi$
(for $1 \le \xi < \omega_1$) which were introduced by D.~Alspach and
S.A.~Argyros~\cite{AA}. Before defining $\S_\xi$ we recall some
general terminology. Let ${\mathcal F}$ be a set of finite subsets of
$\N$. We say that ${\mathcal F}$ is {\em hereditary} if whenever $G
\subseteq F \in {\mathcal F}$ then $G \in {\mathcal F}$.
${\mathcal F}$ is {\em spreading} if whenever $\{ n_1, n_2, \ldots ,
n_k \} \in {\mathcal F}$ with $n_1<n_2< \cdots < n_k$ and $m_1<m_2<
\cdots <m_k$ satisfies $n_i \le m_i$ for $i \le k$ then $\{ m_1, m_2,
\ldots , m_k \} \in {\mathcal F}$.  ${\mathcal F}$ is
{\em pointwise closed} if ${\mathcal F}$ is closed in the topology of
pointwise convergence in $2^{\N}$. ${\mathcal F}$ is called
{\em regular} if it is hereditary, spreading and pointwise closed. If
$A$ and $B$ are two finite subsets of $\N$, then by $A<B$ we mean that
$\max A<\min B$.  Similarly, for $n \in \N$ and $A \subseteq \N$,
$n\le A$ means $n\le\min A$. We assume that $\varnothing < F$ and $F <
\varnothing$ for any non-empty finite set $F\subseteq N$.  If
${\mathcal F}$ and ${\mathcal G}$ are regular then let
$$
{\mathcal F} [ {\mathcal G} ] = \Bigl\{ \bigcup\limits_1^n G_i \mid
n \in \N, \ G_1< \cdots <G_n, \ G_i \in {\mathcal G} \text{ for }i \le n, \
(\min G_i )_1^n \in {\mathcal F} \Bigr\} .
$$
If ${\mathcal F}$ is regular and $n \in \N$ then we define $[
{\mathcal F} ]^n$ by $[{\mathcal F} ]^1= {\mathcal F}$ and $[
{\mathcal F} ]^{n+1} = {\mathcal F} \bigl[ [ {\mathcal F} ]^n \bigr]$.
If $F$ is a finite set then $\# F$ denotes the cardinality of~$F$.  If
$N$ is an infinite subset of $\N$ then $[ N]^{<\infty}$ denotes the
set of all finite subsets of~$N$.  For any ordinal number $1 \le \xi <
\omega_1$, Schreier families $\S _{\xi}$ ($\subseteq [ \N ]^{<
  \infty}$) are defined as follows: set
$$
\S _0=\bigl\{\{n\}\mid n \in \mathbb{N}\bigr\} \cup
\{\varnothing \},
\quad \S_1= \bigl\{ F \subseteq \N\mid \# F \le \min F\bigr\}.
$$
After defining $\S _{\xi}$ for some $\xi < \omega_1 $, set
$$
\S _{\xi +1}= \S_1 [ \S_\xi ] .
$$
If $\xi< \omega_1$ is a limit ordinal and $\S
_\alpha$ has been defined for all $\alpha < \xi$ then fix a sequence
$\xi_n \nearrow \xi$ and define
$$
\S _\xi = \bigl\{ F\mid n \le F \text{ and } F \in S_{\xi_n}\text{ for
    some }n \in \N \bigr\} .
$$
If $N=\{ n_1, n_2, \ldots \}$ is a subsequence of $\N$ with $n_1<n_2<
\cdots$ and ${\mathcal F}$ is a set of finite subsets of $\N$ then we
define ${\mathcal F} (N)=\bigl\{ (n_i)_{i \in F} \mid F \in {\mathcal
  F}\bigr\}$.  We summarize the properties of the Schreier families
that we will need:

\begin{remark} \label{R:Sxi}
\begin{enumerate}
\item\label{reg} Each $\S_\xi$ is a regular family.
\item $\S_\xi \subseteq \S_{\xi +1}$ for every $\xi$. However,
  $\xi<\zeta$ doesn't generally imply $\S_\xi\subseteq \S_\zeta$.
\item\label{monot} Let $1 \le \xi < \zeta < \omega_1$. Then
  there exists $n \in \N$ so that if $n \le F \in \S_\xi$ then $F \in
  \S_\zeta$.
\item\label{prod} For $n,m \in \N$ we have that $\S_n[\S_m]=
  \S_{n+m}$. This fails for infinite ordinals. However, the following
  is true: For all $1 \le \alpha , \beta < \omega_1$ there exist
  subsequences $M$ and $N$ of $\N$ such that $\S_\alpha [ \S_\beta
  ](N) \subseteq \S_{\beta +\alpha}$ and $\S_{\beta + \alpha} (M)
  \subseteq \S_{\alpha}[ \S_\beta]$.  Also for all $1 \le \xi <
  \omega_1$ and $n \in \N$ there exist subsequences $M$ and $N$ of
  $\N$ satisfying $[ \S_\xi ]^n (N) \subseteq \S_{\xi n}$ and
  $\S_{\xi n} (M) \subseteq [ \S_\xi ]^n$.
\item\label{s-seq} Let $1 \le \beta < \alpha < \omega_1$,
  $\varepsilon >0$ and $M$ be a subsequence of $\N$. Then there exists
  a finite set $F \subseteq M$ and $(a_j)_{j \in F} \subseteq
  \mathbb{R}^+$ so that $F \in \S_\alpha (M)$, $\sum_{j \in F} a_j =1$
  and if $G \subseteq F$ with $G \in \S_\beta$ then $\sum_{j \in G}
  a_j < \varepsilon$.
\end{enumerate}
\end{remark}

The proofs can be found in \cite{AGR}.


\section{Classes of strictly singular operators}

Recall that a bounded operator $T$ from a Banach space $X$ to a Banach
 space $Y$ is called \term{strictly singular} if its
restriction to any infinite-dimensional subspace is not an
isomorphism. That is, for every infinite dimensional subspace $Z$ of
$X$ and for every $\varepsilon>0$ there exists $z\in Z$ such that
$\norm{Tz}<\varepsilon\norm{z}$. We say that $T$ is \term{finitely
  strictly singular} if for every $\varepsilon>0$ there exists
$n\in\mathbb N$ such that for every subspace $Z$ of $X$ with $\dim
Z\ge n$ there exists $z\in Z$ such that
$\norm{Tz}<\varepsilon\norm{z}$.  In particular, for $1 \le
p<q\le\infty$ the inclusion operator $i_{p,q}$ from $\ell_p$ to
$\ell_q$ is finitely strictly singular. We will denote by $\K(X,Y)$,
$\SS(X,Y)$, and $\FSS(X,Y)$ the collections of all compact, strictly
singular, and finitely strictly singular operators from $X$ to $Y$,
respectively. If $X=Y$ we will write $\K(X)$, $\SS(X)$, and $\FSS(X)$.
It is known that these sets are norm closed operator ideals in
$\L(X)$, the space of all bounded linear operators on $X$, see
\cite{Milman:70,SSTT} for more details on these classes of operators.
It is well known that $\K(X) \subseteq \FSS(X)\subseteq\SS(X)$.  We
provide the proof for completeness. The second inclusion is obvious.
To prove the first inclusion, suppose that $T$ is not finitely
strictly singular. Then there exists $\varepsilon>0$ and a sequence
$(E_n)$ of subspaces of $X$ such that $\dim E_n=n$ and $T$ satisfies
$\norm{Tx}\ge\varepsilon\norm{x}$ for each $x\in E_n$. Let
$F_n=T(E_n)$. It follows that $\dim F_n=n$ and, for every $n$ and
every $y \in T(S_{E_n})$ we have that $\norm{y} \ge \varepsilon$,
(where $S_{E_n}$ denotes the unit sphere of $E_n$). Let $z_1$ be in
$T(S_{E_1})$. Suppose we have already constructed $z_1,\dots,z_k$ with
$z_i\in T(S_{E_i})$ for $i=1,\dots,k$. Using  \cite[Lemma 1.a.6]{LT1}
or \cite[Lemma of page 2]{Diestel} we can find $z_{k+1}$ in
$T(S_{E_{k+1}})$ such that
$\dist\bigl(z_{k+1},[z_i]_{i=1}^k\bigr)>\frac{\varepsilon}{2}$.
Iterating this procedure we produce a sequence $(z_i)$ in $T(B_X)$
satisfying $\norm{z_i-z_j}>\frac{\varepsilon}{2}$ whenever $i\neq j$.
It follows that $T$ is not compact.

In this article we define and study certain classes of strictly
singular operators. We also refine certain results about strictly
singular operators to the classes of operators that we introduce.

\begin{definition} \label{D:ss} If $X_1$, $X_2$ are Banach spaces, $T
  \in \L (X_1, X_2)$ and $1 \le \xi < \omega_1$, we say that $T$ is
  $\S _\xi$-strictly singular and write $T\in\SS _\xi(X_1, X_2)$ if
  for every $\varepsilon>0$ and every basic sequence $(x_n)$ there
  exist a set $F \in \S _\xi$ and a vector $z\in[x_i]_{i\in F}
  \setminus \{ 0 \}$, ($[x_i]_{i\in F}$ stands for the closed linear
  span of $\{x_i\}_{i\in F}$), such that $\norm{Tz} \le
  \varepsilon\norm{z}$.  If $X_1 =X_2$ then we write $T \in \SS _\xi
  (X_1)$.
\end{definition}

The main difficulty in checking that an operator is $\S_\xi$-strictly
singular, seems to be that one has to verify Definition~\ref{D:ss} for
{\em all} basic sequences $(x_n)$. Notice that without loss of
generality it is enough to check all {\em normalized} basic sequences.
Also notice that if $X_1$, $X_2$ are Banach spaces then $T\in\SS
_\xi(X_1, X_2)$ if and only if for every normalized basic sequence
$(x_n)$ and $\varepsilon >0$ there exist a subsequence $(x_{n_k})$, $F
\in \S _\xi$ and $w \in [x_{n_k}]_{k \in F} \setminus \{ 0 \}$ such
that $\norm{Tw} \le\varepsilon \norm{w}$.  This is easy to see, since
$F \in \S_\xi$ implies that $\{ n_k: k \in F \} \in \S_\xi$.  For
reflexive Banach spaces with bases, we can narrow down even more this
family of basic sequences, as the following remark shows.

\begin{remark}\label{R:SS_xi}
  Let $T \in \L (X_1, X_2)$ and $1 \le \xi < \omega_1$.  If $X_1$ is a
  reflexive Banach space with a basis $(e_n)$ then $T \in \SS_\xi
  (X_1, X_2)$ if and only if for any normalized block sequence $(y_n)$
  of $(e_n)$ and $\varepsilon >0$ there exists $G \in \S _\xi$ and $w
  \in [y_n]_{n \in G} \setminus \{ 0 \}$ such that $\norm{Tw} \le
  \varepsilon \norm{w}$.
\end{remark}

Remark~\ref{R:SS_xi} follows from the following classical fact.  Two
basic sequences $(x_n)$ and $(y_n)$ are called $C$-equivalent for some
$C>1$, denoted by $(x_n) \eqv{C} (y_n)$, if for every $(a_n) \in
c_{00}$ we have that $\norm{\sum{a_n}x_n} \eqv{C}\norm{\sum{a_n}y_n}$.
(We write $a \eqv{C}b$ if $\frac{1}{C} a \le b \le C a$.)  Two basic
sequences $(x_n)$ and $(y_n)$ are called equivalent, denoted by $(x_n)
\approx (y_n)$, if they are $C$-equivalent for some $C \ge 1$.  Since
$X_1$ is reflexive then every normalized basic sequence $(x_n)$ in
$X_1$ is weakly null and therefore by~\cite{BP} it has a subsequence
$(x_{n_k})$ which is equivalent to a block sequence $(y_k)$ of $(e_n)$
and $\norm{x_{n_k}-y_k}\to 0$.


\begin{remark}\label{subseq}
  Let $(x_n)$ be a bounded sequence in a Banach space $X$. Then there
  is a subsequence $(x_{n_k})$ such that one of the following
  conditions hold.
  \begin{enumerate}
  \item\label{conv} $(x_{n_k})$ converges;
  \item\label{l1} $(x_{n_k})$ is equivalent to the unit vector basis of
    $\ell_1$;
  \item\label{dbs} The difference sequence $(d_k)$ defined by
    $d_k=x_{n_{2k+1}}-x_{n_{2k}}$ is a seminormalized weakly null
    basic subsequence.  Moreover, if $X$ has a basis then $(d_k)$ is
    equivalent to a block sequence of the basis.
  \end{enumerate}
\end{remark}

This is a standard result. Indeed, if $(x_n)$ has no subsequences
satisfying~\eqref{conv} or~\eqref{l1} then Rosenthal's
$\ell_1$~Theorem yields a weakly Cauchy subsequence $(x_{n_k})$. By
passing to a further subsequence we may assume that the sequence
$(x_{n_{k+1}}-x_{n_k})$ is weakly null and seminormalized.
Now~\eqref{dbs} follows by~\cite{BP}.

In view of this Remark~\ref{subseq}, the requirement ``every basic
sequence'' in Definition~\ref{D:ss} is ``almost'' as general as
``every sequence''.

\begin{proposition} \label{P:SS_xi} Suppose that $X$ and $Y$ are two
  Banach spaces and $1\le\xi , \zeta <\omega_1$. Then
  \begin{enumerate}
  \item\label{SSa} $\FSS(X,Y)\subseteq\SS_\xi(X,Y)\subseteq\SS(X,Y)$.
  \item\label{SSb} If $1 \le \xi < \zeta < \omega_1$ then $\SS_\xi(X,Y)
    \subseteq \SS_\zeta(X,Y)$.
  \item\label{SSnorm} $\SS _\xi(X)$ is norm-closed;
  \item\label{SSprod} If $S\in\SS_\xi(X)$ and $T\in\L(X)$ then
    $TS$ and $ST$ belong to $\SS_\xi(X)$.
  \item\label{SSsum} If $S \in\SS_\xi(X)$ and $T \in \SS_\zeta (X)$
    then $S+T\in\SS_{\xi + \zeta}(X)$. In particular, if $S,T \in
    \SS_\xi(X)$ then $S+T\in\SS_{\xi 2}(X)$.
  \end{enumerate}
\end{proposition}

\begin{proof} (\ref{SSa}) It is obvious.

  (\ref{SSb}) Indeed, for $1 \le \xi < \zeta < \omega_1$, by
  Remark~\ref{R:Sxi}\eqref{monot} there exists $N \in \N$ such that if
  $\S_\xi \cap \bigl[ \{ N, N+1, \ldots \} \bigr]^{< \infty} \subseteq
  \S_\zeta$.  Now if $T \in \SS_\xi(X,Y)$, $\varepsilon$ is a positive
  number and $(x_n)$ is a normalized basic sequence in $X$ then
  consider the basic sequence $(y_n)$ where $y_i=x_{N+ i}$. There
  exists $F \in \S_\xi$ and $z \in [y_i]_{i \in F} \setminus \{ 0 \} $
  such that $\norm{Tz} \le \varepsilon \norm{z}$. Since $F \in \S_\xi$
  and $F \subseteq \{ N, N+1, \ldots \}$ we have that $F \in
  \S_\zeta$.

  (\ref{SSnorm}) Let $(T_n)_n \subset \SS _\xi(X)$, $T \in \L (X)$ and
  $\lim_n T_n = T$. Let $(x_n)$ be a seminormalized basic sequence in
  $X$, and $\varepsilon >0$. Let $n_0 \in \N$ such that
  $\norm{T_{n_0} - T } \le \varepsilon /2$.  Since $T_{n_0} \in \SS
  _\xi(X)$, there exists $F \in \S _\xi$ and $z \in [x_i]_{i \in F}
  \setminus \{ 0 \} $ such that $\norm{T_{n_0} z}\le
  \frac{\varepsilon}{2} \norm{z}$. Thus
  $$
  \norm{Tz} \le \bignorm{(T_{n_0}-T)z} + \norm{T_{n_0} z} \le
  \tfrac{\varepsilon}{2} \norm{z} + \tfrac{\varepsilon}{2} \norm{z} =
  \varepsilon \norm{z}.
  $$

  (\ref{SSprod}) Let $S\in\SS_\xi(X)$ and $T\in\L(X)$. We show that
  $TS \in \SS_\xi(X)$. Let $(x_n)$ be a basic sequence in $X$ and
  $\varepsilon>0$. If $T=0$ then it is obvious that $TS\in\SS_\xi
  (X)$.  Suppose that $T\neq 0$, then there exists $F\in\S_\xi$ and
  $z\in[x_n]_{n\in F} \setminus \{ 0 \} $ such that
  $\norm{Sz}\le\frac{\varepsilon}{\norm{T}}\norm{z}$. Thus,
  $\norm{TSz}\le\norm{T}\norm{Sz}\le\varepsilon\norm{z}$.  The proof
  that $ST \in \SS_\xi(X)$ is due to A.~Popov \cite{Popov} who
  improved our original argument which only worked in reflexive
  spaces.

  (\ref{SSsum}) Let $(x_n)_{n \in \N}$ be a normalized basic sequence
  and $\varepsilon>0$.  By Remark~\ref{R:Sxi}\eqref{prod} let
  $N=(n_i)$ be a subsequence of $\N$ such that $\S_\zeta [\S_\xi ] (N)
  \subseteq \S_{\xi + \zeta}$.  Find $F_1\in\S_\xi$ and
  $w_1\in[x_{n_i}]_{i\in F_1}$ such that $\norm{w_1}=1$ and
  $\norm{Sw_1}<\frac{\varepsilon}{8C}$ where $C$ is the basis constant
  of $(x_n)$.  Since $(x_{n_i})_{i>F_1}$ is again a basic sequence, we
  can find $F_2\in\S_\xi$ and $w_2\in[x_{n_i}]_{i\in F_2}$ such that
  $F_1<F_2$, $\norm{w_2}=1$, and
  $\norm{Sw_2}<\frac{\varepsilon}{16C}$. Proceeding inductively we
  produce sets $F_1<F_2<\dots$ and vectors
  $w_k\in[x_{n_i}]_{i\in F_k}$ with $\norm{w_k}=1$ and
  $\norm{Sw_k}<\frac{\varepsilon}{2^{k+2}C}$.  Since $(w_k)$ is a
  basic sequence, we find $G\in\S_\zeta$ and $z\in[w_k]_{k\in G}
  \setminus \{ 0 \} $ such that $\norm{Tz}\le \frac{\varepsilon}{2}
  \norm{z}$.  Suppose that $G=\{k_1,\dots,k_m\}$ and
  $z=\sum_{i=1}^ma_iw_{k_i}$. Then we can write
  $z=\sum_{i\in F}b_ix_{n_i}$ for some $F\in \S_\xi [\S_\xi]$. By the
  choice of $N$ we have that $z \in [(x_i)_{i \in H}]$ for some $H \in
  \S_{\xi + \zeta}$. Also, $\abs{a_i}\le 2C\norm{z}$. It follows that
  $\norm{Sz}\le\sum_{i=1}^m\abs{a_i}\norm{Sw_{k_i}}\le
  2C\frac{\varepsilon}{4C} \norm{z}= \frac{\varepsilon}{2} \norm{z}$,
  so that $\norm{(S+T)z} \le \varepsilon\norm{z}$.
\end{proof}

Of course, if $1 \le p< q < \infty$ then any bounded operator from
$\ell_q$ to $\ell_p$ is compact.  Also every bounded operator from
$\ell_p$ to $\ell_q$ is strictly singular, \cite{LT1}.

\begin{example} \label{pq}
  \textit{If $1 \le p< q < \infty$ then any bounded operator $T \in \L
    (\ell_p, \ell_q)$ belongs to $\SS_1 (\ell_p, \ell_q)$.}

  If $1 < p$ then we can apply Remark~\ref{R:SS_xi}. Let $(x_n)$ be a
  normalized block sequence in $\ell_p$ and $\varepsilon >0$.  If
  $\inf_i \norm{Tx_i}_q =0$ then we are done (we denote by
  $\norm{ \cdot}_p$ and $\norm{\cdot}_q$ the norms of $\ell_p$ and
  $\ell_q$ respectively), hence assume that $(Tx_n)$ is
  seminormalized.  Since $(x_n)$ is weakly null, $(Tx_n)$ is weakly
  null. By standard gliding hump arguments \cite{BP} we can pass to a
  subsequence $(Tx_{n_i})$ such that for some seminormalized block
  sequence $(y_n)$ in~$\ell_q$,
  $$
    \Bignorm{ \sum_i a_i Tx_{n_i}}_q \le 2 \Bignorm{ \sum_i a_i y_i}_q
     \text{ for every }(a_i) \in c_{00}.
  $$
  Hence for $\varepsilon >0$ one can choose $N \in \N$ such that
  $$
  \Bignorm{ T \Bigl( \sum_{i=1}^N x_{n_{N+i}} \Bigr)}_q \le
  \varepsilon \Bignorm{\sum_{i=1}^N x_{n_{N+i}}}_p.
  $$
  Suppose that $p=1$. Let $(x_n)$ be a normalized basic sequence in
  $\ell_1$ and $\varepsilon >0$. By H.P.~Rosenthal's $\ell_1$ theorem
  \cite{R} after passing to a subsequence and relabeling we can assume
  that $(x_n)$ is $K$-equivalent to the unit vector basis of $\ell_1$
  for some $K < \infty$.

  By applying Remark~\ref{subseq} to $(Tx_n)$ there exists a
  subsequence $(x_{n_k})$ of $(x_n)$ such that the sequence $(d_k)$
  defined by $d_k=Tx_{n_{2k+1}}-Tx_{n_{2k}}$ is either norm null or
  satisfies (\ref{dbs}) of Remark~\ref{subseq}. If $(d_k)$ is norm
  null, then there exists $m\ge 2$ such that
  $\norm{d_m}<\frac{2\varepsilon}{K}$, so that
  $$
    \norm{Tx_{n_{2m+1}}-Tx_{n_{2m}}}<
     \varepsilon\cdot\tfrac{2}{K}\le
     \varepsilon\norm{x_{n_{2m+1}}-x_{n_{2m}}}.
  $$
  Since $\{ n_{2m},n_{2m+1} \}\in\S_1$ we have $T \in \SS_1(\ell_1, \ell_q)$.

  If $(d_k)$ is $C$-equivalent to a block sequence of the standard
  basis of $\ell_q$ then for $\varepsilon >0$ one can choose $N\in\N$
  such that
  $$
  \Bignorm{T\left( \sum_{k=1}^N(x_{n_{2(N +k)+1}}- x_{n_{2(N+k)}})
    \right) }
  \le \varepsilon
  \Bignorm{\sum_{k=1}^N(x_{n_{2(N +k)+1}}- x_{n_{2(N+k)}})}.
  $$
\end{example}

\begin{example}
  Suppose $1<p<q<\infty$ with $p\neq q$. Then it is known (see
  \cite{Milman:70,P,SSTT}) that
  $\FSS(\ell_p,\ell_q)\neq\SS(\ell_p,\ell_q)=\L(\ell_p,\ell_q)$.
  Therefore, Example~\ref{pq} yields
  $\FSS(\ell_p,\ell_q)\neq\SS_1(\ell_p,\ell_q)$.
\end{example}

\begin{example} \label{E:SS_xi}
  \textit{An example of a space $X$ where $\SS_\xi(X)\neq\SS_\zeta(X)$
    for some $1\le\xi<\zeta<\omega_1$.}

  Fix $1 \le \xi < \omega_1$ consider the
  space $T[\S_\xi , \frac{1}{2}]$ which is the completion of $c_{00}$
  with the norm that satisfies the implicit equation:
  $$
    \norm{x}_\xi =
    \max\Bigl\{\norm{x}_\infty,\sup\tfrac{1}{2}\sum_{i}\norm{E_ix}_\xi\Bigr\},
  $$
  where $\norm{ \cdot }_\infty$ stands for the $\ell_\infty$ norm, and
  the supremum is taken for all sets $E_1< E_2< \cdots $ such that
  $(\min E_i )_i \in \S_\xi$.

  Since $\xi w$ is a limit ordinal, without loss of generality we can
  assume that the sequence of
  ordinals in the definition of $\S_{\xi w}$ starts with
  $\xi$, then $\S_\xi\subseteq\S_{\xi w}$ and, therefore,
  $T[\S_{\xi \omega},\frac{1}{2}]\subseteq T[\S_\xi , \frac{1}{2}]$.
  Consider the inclusion operator $i_\xi: T[\S_{\xi \omega},
  \frac{1}{2}] \to T[\S_\xi , \frac{1}{2}]$.  Then $i_\xi \in
  \SS_{\xi \omega}\bigl(T[\S_{\xi \omega}, \frac{1}{2}], T[\S_\xi ,
  \frac{1}{2}]\bigr)$ but $i_\xi \not \in \SS_\xi \bigl(T[\S_{\xi \omega},
  \frac{1}{2}], T[\S_\xi , \frac{1}{2}]\bigr)$.

  Indeed, it is easy to verify that $i_\xi \not \in \SS_\xi
  \bigl(T[\S_{\xi \omega}, \frac{1}{2}], T[\S_\xi ,
  \frac{1}{2}]\bigr)$, since for every $F \in \S_\xi$ and scalars
  $(a_i)_{i \in F}$, we have that
  $$
  \Bignorm{i_\xi\bigl(\sum_{i \in F} a_i e_i\bigr) }_\xi =
    \max \Bigl\{ \max_{i \in F} |a_i|, \tfrac{1}{2}\sum_{i \in F}\abs{a_i}\Bigr\} =
  \bignorm{\sum_{i \in F} a_i e_i}_{\xi \omega},
  $$
  where $(e_i)$ denotes the standard basis of $T[\S_{\xi \omega}, \frac{1}{2}]$.

  Now we verify that
  $i_\xi \in \SS_{\xi \omega}(T[\S_{\xi \omega}, \frac{1}{2}],
  T[\S_\xi , \frac{1}{2}])$.  First recall that
  $T[\S_{\xi \omega}, \frac{1}{2}]$ is a reflexive Banach space with a
  basis \cite[Proposition~1.1]{AD}. Thus we can apply
  Remark~\ref{R:SS_xi}. Let $(x_n)$ be a normalized block
  sequence in $T[\S_{\xi \omega}, \frac{1}{2}]$ and $\varepsilon >0$.
  If there exists $n \in \N$ such that
  $\norm{i_\xi x_n }_\xi = \norm{x_n}_\xi \le \varepsilon,$
  then we are done. Else assume that $(i_\xi x_n)_n$ is
  seminormalized.  Let $n_i = \min \supp x_i$ (with respect to $(e_i)$). By
  \cite[Proposition~4.10]{LM} we have that $(i_ \xi x_i) \eqv{C}
  (e_{n_i})$ where $C := 96 \sup_i \norm{x_i}_\xi / \inf_i
  \norm{x_i}_\xi$.  We have the following claim which uses the idea
  and generalizes \cite[Proposition~1.5]{AD}.

  \textit{Claim~1:} For every $\eta >0$ there exists $F \in
  \S_{\xi \omega}$ and a convex combination $x :=\sum_{i \in F}
  a_{n_i} e_{n_i}$ such that $\norm{x}_\xi < \eta$.

  Once Claim~1 is proved then by letting
  $\eta:=\frac{\varepsilon}{2C}$ it follows that
  \begin{displaymath}
    \Bignorm{ i_\xi \bigl( \sum_{i \in F} a_{n_i} x_i \bigr) }_\xi =
    \Bignorm{ \sum_{i \in F} a_{n_i} x_i }_\xi \le
    C \norm{ x}_\xi\le \tfrac{\varepsilon}{2} =
    \tfrac{\varepsilon}{2} \sum_{i \in F} a_{n_i} \norm{x_i}_{\xi \omega}
    \le \varepsilon \Bignorm{ \sum_{i \in F} a_{n_i} x_i}_{\xi \omega}.
  \end{displaymath}

  Thus it only remains to establish Claim~1. For this purpose we need
  to identity a norming set $N^\xi$ of
  $T[\S_\xi , \frac{1}{2}]$. We follow \cite[page~976]{AD}: Let
  $$
    N^\xi_0= \{ \pm e_n^* : n \in \N \}\cup\{0\}.
  $$
  If $N^\xi_s$ has been defined for some $s \in \N \cup \{ 0 \}$, then
  we define
  \begin{multline*}
    N^\xi_{s+1}  = N^\xi_s \cup \Bigl\{ \tfrac{1}{2} (f_1 + \cdots +f_d):
        f_i \in N^\xi_s,\, (i=1, \ldots , d),\\
     \supp f_1 < \supp f_2 < \cdots < \supp
    f_d \text{ and } (\min\supp f_i)_{i=1}^d \in \S_\xi \Bigr\} .
  \end{multline*}
  Finally set $N^\xi = \cup_{s=0}^\infty N^\xi_s$ and the set $N^\xi$
  is a norming set for $T[\S_\xi , \frac{1}{2}]$, i.e. we have
  $\norm{x}_\xi = \sup_{ x^* \in N^\xi} x^*(y) $ for all $x \in
  T[\S_\xi , \frac{1}{2}]$.

  Now we prove Claim~1. First choose $\ell \in \N$ such that
  $\frac{1}{2^ \ell} < \frac{\eta}{2}$.  We have the following claim
  which follows immediately from Remark~\ref{R:Sxi}\eqref{s-seq}.

 \textit{Claim 2:} There exists a convex combination $x
  =\sum_{i \in F} a_{n_i} e_{n_i}$ such that $F \in \S_{\xi \ell+1} \cap
  \S_{\xi \omega}$ and $\sum_{i \in G} a_i < \frac{\eta}{2}$ for all
  $G \in \S_{\xi \ell}$.

  Let $x$ as in Claim~2. In order to estimate
  $\norm{x}_\xi$ from above, let $x^* \in N^\xi$. Let $L := \{ k \in
  \N : \abs{x^*(e_k)} \ge \frac{1}{2^\ell} \}$. Then $L \in
  \S_{\xi \ell}$.  Therefore
  \begin{displaymath}
    \bigabs{ x^* (x)} \le \bigabs{(x^*|_L)(x)}+ \bigabs{(x^*|_{L^c})(x)}
     \le \sum_{k \in L} a_k + \frac{1}{2^\ell} < \tfrac{\eta}{2} +
    \tfrac{\eta}{2} = \eta.
  \end{displaymath}
  This finishes the proof of Claim~1 and the proof that $i_\xi \in
  \SS_{\xi \omega}(T[\S_{\xi \omega}, \frac{1}{2}], T[\S_\xi ,
  \frac{1}{2}])$.
\end{example}

\begin{remark}
  Suppose that $X$ and $Y$ are Banach spaces, $1 \le \xi < \omega_1$ and $T \in
  \SS_\xi (X,Y)$. Let $\widetilde T\in\L(X\oplus Y)$ given by
  $(x,y)\mapsto(0,Tx)$, that is,
  \begin{math}
    \widetilde T=\left(
    \begin{smallmatrix}
     0 & 0 \\ T & 0
    \end{smallmatrix}\right).
  \end{math}
  Then $\widetilde T\in\SS_\xi(X\oplus Y)$. Conversely, if $\widetilde
  T \in \SS_\xi(X\oplus Y)$ then $T \in \SS_\xi (X,Y)$.
\end{remark}

The converse is obvious. To see the forward implication, pick a
normalized basic sequence $(x_n,y_n)$ in $X\oplus Y$ and
$\varepsilon>0$. Since $(x_n)$ is bounded, there exists a subsequence
$(x_{n_i})$ of $(x_i)$ which satisfies one of the options in
Remark~\ref{subseq}. Set $d_k=x_{n_{2k+1}}-x_{n_{2k}}$.

In case (\ref{conv}), $d_m\to 0$, so we can choose $m$ such that
$\norm{d_m}<\frac{\varepsilon}{C\norm{T}}$, where $C$ is the basis
constant of $(x_n,y_n)$. Put
$h=(x_{n_{2m+1}},y_{n_{2m+1}})-(x_{n_{2m}},y_{n_{2m}})$, then $\supp
h=\{ n_{2m}, n_{2m+1} \} \in\S_\xi$ and
$$\norm{\widetilde{T}h}=\bignorm{(0,T(x_{n_{2m+1}}-x_{n_{2m}}))}\le\norm{T}\norm{d_m}
<\frac{\varepsilon}{C}\le\varepsilon\norm{h}.$$

In case (\ref{l1}), since $T \in \SS_\xi (X,Y)$ and $(x_n)$ is a basic
sequence, we can find $F\in\S_\xi$ and non-zero scalars
$(a_i)_{i\in F}$ such that if $w=\sum_{i\in F}a_ix_{n_i}$ then
$\norm{Tw}\le\varepsilon\norm{w}$. Let
$h=\sum_{i\in F}a_i(x_{n_i},y_{n_i})$, then
\begin{displaymath}
   \norm{\widetilde{T}h}=\bignorm{(0,Tw)}\le\varepsilon\norm{w}\le\varepsilon\norm{h},
\end{displaymath}
where, without loss of generality, we assume that
$\bignorm{(0,y)}=\norm{y}$ for all $y \in Y$.

In case (\ref{dbs}), suppose that $(d_k)$ is a basic seminormalized
sequence. Then there exists $G\in\S_\xi$ and non-zero scalars
$(a_k)_{k\in G}$ such that $\norm{Tw}\le\varepsilon\norm{w}$ where
$w=\sum_{k\in G}a_kd_k$. Set
$$h=\sum_{k\in G}a_k\bigl((x_{n_{2k+1}},y_{n_{2k+1}})-(x_{n_{2k}},y_{n_{2k}})\bigr).$$
Then
$\norm{\widetilde{T}h}=\bignorm{(0,Tw)}\le\varepsilon\norm{w}\le\varepsilon\norm{h}$.
It is left to show that $\supp h\in\S_\xi$.  For a set $A \subseteq \N$ define
$A^{\times 2}= \cup_{i \in A} \{ 2 i, 2i +1 \}$. By transfinite induction it is easy to see
that if $A \in \S_\xi$ then $A^{\times 2} \in \S_\xi$. Thus $F:= G^{\times 2} \in \S_\xi$.
Therefore, $\supp h= \{ n_k \mid  k \in F \} \in\S_\xi$ since $\S_\xi$ is spreading
Remark~\ref{R:Sxi}\eqref{reg}.

\section{Schreier-spreading sequences and some equivalence relations}

Recall the notion of \term{spreading model}.  It is  shown in
\cite{BS1,BS2} that for every seminormalized basic sequence $(y_i)$ in
a Banach space and for every $\varepsilon_n \searrow 0$ there exists a
subsequence $(x_i)$ of $(y_i)$ and a seminormalized basic sequence
$(\tilde{x}_i)$ (in another Banach space) such that for all $n\in
\mathbb{N}$, $(a_i)^n_{i=1}\in [-1,1]^n$ and $n\le k_1 <\ldots <k_n$
one has
\begin{equation}\label{sm}
  \Bigabs{ \bignorm{ \sum\limits_{i=1}^n a_ix_{k_i} } -
    \bignorm{ \sum\limits_{i=1}^n a_i\tilde x_i} } <\varepsilon_n.
\end{equation}
The sequence $(\tilde x_i)$ is called the
\term{spreading model of $(x_i)$} and it is a suppression-1
unconditional basic sequence if $(y_i)$ is weakly null. We refer the
reader to \cite{BS1}, \cite{BS2} and \cite[I.3.~Proposition 2]{BL} for
more information about spreading models. Spreading models of weakly
null seminormalized basic sequences have been studied in \cite{AOST},
where for a Banach space~$X$, the set of all spreading models of all
seminormalized weakly null basic sequences of $X$ is denoted by $\SP
_w(X)$.  Also $\#\SP_w (X)$ denotes the cardinality of the quotient of
$\SP_w (X)$ with respect to the equivalence relation $\approx$. In
other words, $\#\SP_w (X)$ is the largest number of pair-wise
non-equivalent spreading models of weakly null seminormalized basic
sequences in $X$,~(\cite{AOST}).

We will use the following standard fact whose proof is left to the reader.

\begin{lemma}\label{sm-spread}
  Suppose that $(x_n)$ is a seminormalized basic sequence with a
  spreading model $(\tilde x_n)$. Then, for every $\varepsilon >0$
  there exists $n_0\in\mathbb N$ such that
  \begin{displaymath}
    \Bignorm{\sum_{i=1}^na_i\tilde x_i}\eqv{1 + \varepsilon}
    \Bignorm{\sum_{i=1}^na_i x_{k_i}}
  \end{displaymath}
  whenever $n_0\le n\le k_1<\dots<k_n$ and $a_1,\dots,a_n\in\mathbb R$.
\end{lemma}

Motivated by the definition of spreading model we now define
the notion of a \term{Schreier spreading sequence}.

\begin{definition}
  Let $X$ be a Banach space. We say that a seminormalized basic
  sequence $(x_n)$ in $X$ is \term{Schreier spreading}, if there
  exists $1 \le C < \infty$ such that for every $F =\{ f_1, f_2,
  \ldots , f_n \} , G=\{ g_1, g_2, \ldots , g_n \} \in \S_1$ and
  scalars $(a_i)_{i =1}^n$ we have
\begin{displaymath}
  \bignorm{ \sum\limits_{i =1}^n a_i x_{f_i} } \eqv{C}
  \bignorm{ \sum\limits_{i =1}^n a_i x_{g_i} }.
\end{displaymath}
Let $\SP _{1,w} (X)$ denote the set of seminormalized weakly null
basic sequences in $X$ which are Schreier spreading (here
the index ``$1$'' reminds us of $\S_1$, and the index ``$w$'' reminds us of
\textit{weakly null}).
\end{definition}

It follows immediately from the results of Brunel and Sucheston \cite{BS1,BS2}
and  Lemma~\ref{sm-spread} that

\begin{remark}\label{BS}
  Every seminormalized basic sequence has a Schreier spreading
  subsequence.
\end{remark}

Now for $1 \le \xi < \omega_1$ we define equivalence relations
$\approx_\xi$ on $\SP_{1,w}(X)$ as follows:

\begin{definition}
  Let $X$ be a Banach space and $1 \le \xi < \omega_1$. Define an
  equivalence relation $\approx_\xi$ on $\SP_{1,w}(X)$ as follows:
  if $(x_n)$ and $(y_n)$ are two Schreier spreading sequences in $X$,
  we write $(x_n) \approx_\xi (y_n)$ if there exists $1 \le K <
  \infty$ such that for every $F \in \S_\xi$ and scalars
  $(a_i)_{i \in F}$ we have that
  $$
    \Bignorm{\sum_{i \in F} a_i x_i } \eqv{K} \Bignorm{ \sum_{i \in F} a_i y_i}.
  $$
\end{definition}

\begin{proposition} \label{approx}
  Suppose that $(x_n)$ is a Schreier spreading seminormalized basic
  sequence in $X$.
  \begin{enumerate}
  \item\label{app-a} If $(x_{n_k})$ is a subsequence of $(x_n)$ then
    $(x_{n_k})_k$ is Schreier spreading and $(x_n) \approx_1
    (x_{n_k})$.
  \item\label{app-y} If $(x_n)\approx_1(y_n)$ for another basic
    sequence $(y_n)$, then $(y_n)$ is Schreier spreading;
  \item\label{app-b} There exists a normalized Schreier spreading
    sequence $(y_n)$ in $X$ such that $( x_n) \approx (y_n)$.
  \item\label{app-c} If $X$ is a reflexive space with a basis $(e_n)$
    then there exists a seminormalized block sequence $(y_n)$ of
    $(e_n)$ such that $(y_n)$ is Schreier spreading and $( x_n)
    \approx_1 (y_n)$.
  \end{enumerate}
\end{proposition}

\begin{proof}
  (\ref{app-a}) and (\ref{app-y}) are trivial.

  (\ref{app-b}) By standard perturbation arguments, one can find
  $c_0\in\bigl[\inf_n\norm{x_n},\sup_n\norm{x_n}\bigr]$ and a subsequence
  $(x_{n_k})$ of $(x_n)$ such that
  $$
    \Bigl(\frac{x_{n_k}}{\norm{x_{n_k}}}\Bigr) \approx \Bigl(
    \frac{x_{n_k}}{c_0}\Bigr).
  $$
Hence, if $y_k=\frac{x_{n_k}}{\norm{x_{n_k}}}$ then $(y_k)$ is
  normalized basic Schreier spreading and $( y_k) \approx ( x_{n_k})
  \approx_1 (x_k)$.

  (\ref{app-c}) Since $X$ is reflexive, $(x_n)$ is weakly
  null. A standard gliding hump argument yields a
  subsequence $(x_{n_k})$ of $(x_k)$ and a block sequence $(y_k)$ of
  $(e_k)$ such that $(x_{n_k}) \approx (y_k)$ which obviously
  implies the result since $(x_k) \approx_1 (x_{n_k})$.
\end{proof}

\begin{corollary}
  For every Banach space $X$ we have
  $$\#\SP_w(X)=\#\bigl(\SP_{1,w}(X)/\!\!\approx_1\bigr)\le\nSP.$$
\end{corollary}

\begin{proof}
  The inequality $\#\bigl(\SP_{1,w}(X)/\!\!\approx_1\bigr)\le\nSP$ is
  obvious. To show that
  $\#\SP_w(X)=\#\bigl(\SP_{1,w}(X)/\!\!\approx_1\bigr)$ we define a
  bijection $\Phi$ from the set of $\approx$-equivalence classes of
  $\SP_w(X)$ to the set of $\approx_1$-equivalence classes of
  $\SP_{1,w}(X)$. Suppose that $(\tilde x_n)\in\SP_w(X)$ is the
  spreading model of a weakly null seminormalized basic sequence
  $(x_n)$. Then by Lemma~\ref{sm-spread} there exists $n_0 \in \N$
  such that $(x_n)_{n \ge n_0} \in \SP_{1,w}(X)$ and $(x_n)_{n \ge
    n_0} \approx_1 (\tilde{x}_n)_{n \in \N}$. Define $\Phi\colon\bigl((\tilde
  x_n)/\!\!\approx\bigr)\mapsto\bigl((x_n)_{n \ge
    n_0}/\!\!\approx_1\bigr)$.  Obviously $\Phi$ is well defined and
  one-to-one.  It follows from Remark~\ref{BS} and
  Proposition~\ref{approx}(\ref{app-a}) that $\Phi$ is onto.
\end{proof}

\section{Compact products}

Milman \cite{Milman:70} proved that the product of any two strictly
singular operators in $L_p [0,1]$ ($1 \le p < \infty$) or $C[0,1]$ is
compact. In this section we extend the techniques used by Milman to
spaces with finite $\nSPxi$.

\begin{theorem}\label{main}
  Let $X$ be a Banach space, $1 \le \xi < \omega_1$ and $n \in \N \cup
  \{ 0 \}$.  If $\nSPxi=n$, $S\in\SS(X)$,
  and $T_1,\dots,T_n \in\SS_\xi (X)$, then $T_nT_{n-1}\dots T_1S$ is
  compact.  Moreover, if $\ell_1$ does not isomorphically embed in $X$
  then $T_nT_{n-1}\dots T_1$ is compact.

  Furthermore, if $\nSP=n$, and $T_1,\ldots,T_{n+1}\in\SS(X)$, then
  $T_{n+1}T_n\dots T_1$ is compact.  Moreover, if $\ell_1$ does not
  isomorphically embed in $X$ then $T_nT_{n-1}\dots T_1$ is compact.
\end{theorem}

\begin{proof}
  For simplicity, we present the proof in the case $n=2$. However, it
  should be clear to the reader how to extend the proof to $n>2$ or
  $n=1$. The case $n=0$ will be treated at the end.
  Thus, for the sake of contradiction, suppose that the
  conclusion of the theorem fails, i.e., $T_2T_1S$ is not compact or
  $\ell_1\not\hookrightarrow X$ and $T_2T_1$ is not compact.

  \textit{Claim:} There exists a seminormalized weakly Cauchy
  sequence $(u_n)$ such that $(T_2T_1u_n)$ has no convergent
  subsequences.

  If $\ell_1\not\hookrightarrow X$ and $T_2T_1$ is not compact then
  one can find a normalized sequence $(u_n)$ in $X$ such that
  $(T_2T_1u_n)$ has no convergent subsequences. By Rosenthal's Theorem
  \cite{R} we can assume that $(u_n)$ is
  weakly Cauchy.

  Suppose now that $T_2T_1S$ is not compact. Again, find a normalized
  sequence $(v_n)$ in $X$ such that $(T_2T_1Sv_n)$ has no convergent
  subsequences. Put $u_n=Sv_n$. Note that $(T_2T_1u_n)$ has no
  convergent subsequences, so that $(u_n)$ is seminormalized.  Apply
  Rosenthal's Theorem to $(v_n)$. If $(v_n)$ has a weakly Cauchy
  subsequence then, by passing to this subsequence, $(u_n)$ is also
  weakly Cauchy, and we are done. Suppose not, then, by passing to a
  subsequence and relabeling, we can assume that $(v_n)$ is equivalent
  to the unit vector basis of $\ell_1$. Now apply Rosenthal's Theorem
  to $(u_n)$.  If $(u_n)$ has a subsequence equivalent to the unit
  vector basis of $\ell_1$ then, after passing to this subsequence and
  relabeling, we would get that the restriction of $S$ to
  $[v_n]_{n=1}^\infty$ is equivalent to an isomorphism on $\ell_1$,
  which contradicts $S$ being strictly singular. Therefore, $(u_n)$
  must have a weakly Cauchy subsequence.  This completes the proof of
  the claim.

  Since $(T_2T_1u_n)$ has no convergent subsequences, by passing to a
  subsequence and relabeling, we can assume that $(T_2T_1u_n)$ is
  $\varepsilon$-separated for some $\varepsilon >0$.
  Thus the sequences $(x_n)$, $(y_n)$ and $(z_n)$ are seminormalized,
  where $x_n:=u_{n+1}-u_n$, $y_n:=T_1x_n$ and $z_n:=T_2T_1x_n$.
  Since $(u_n)$ is weakly Cauchy, it follows that $(x_n)$, $(y_n)$, and
  $(z_n)$ are weakly null. By using Corollary~1 of~\cite{BP} and
  Remark~\ref{BS}, pass to subsequences and relabel in order to assume
  that $(x_n)$, $(y_n)$, and $(z_n)$ are basic and Schreier spreading.

  Since $T_1 (x_n)=y_n$ for all $n$ and $T_1 \in \SS_\xi (X)$ we have
  that $(x_n) \not \approx_\xi (y_n)$. Similarly, since $T_2(y_n)=z_n$
  for all $n$ and $T_2 \in \SS_\xi(X)$ we obtain that $(y_n) \not
  \approx_\xi (z_n)$. Finally by Proposition~\ref{P:SS_xi}(\ref{SSprod}),
  we have that $T_2T_1 \in \SS_\xi (X)$, so that $(x_n) \not
  \approx_\xi (z_n)$. Thus $\nSPxi\ge 3$, which is a contradiction.

  For the ``furthermore'' statement, if $\nSP=2$ then we can modify
  the above proof to merely assume that $T_1,T_2 \in \SS (X)$. Notice
  that since $T_1 (x_n)=y_n$ for all $n$ and $T_1 \in \SS (X)$ we have
  that $(x_n) \not \approx (y_n)$.  Indeed, otherwise $T_1$ induces
  the restriction operator from $[(x_n)]$ to $[(y_n)]$ via
  $\sum_{i=1}^\infty a_nx_n\mapsto\sum_{i=1}^\infty a_ny_n$. This
  restriction is one-to one since $(y_n)$ is a basic sequence, and
  onto since $(x_n)\approx (y_n)$. Hence, the restriction of $T$ to
  $[x_n]$ would be an isomorphism, contradiction.  Similarly, $(y_n)
  \not \approx (z_n)$ and $(x_n) \not \approx (z_n)$. Thus $\nSP\ge 3$
  which is a contradiction.

  The statement as well as the proof of this result for $n=0$ should
  be given special attention. The assumptions $\nSPxi=0$ or $\nSP=0$,
  combined with Remark~\ref{BS}, simply mean that there is no
  seminormalized weakly null basic sequence in $X$.  The conclusion of
  the statement if $n=0$ simply means that $\K(X)= \SS (X)$. In order
  to verify the result, if $S \in \SS (X) \backslash \K (X)$ then
  there exists a normalized sequence $(v_n)$ such that $(S v_n)$ has
  no convergent subsequence. By Remark~\ref{subseq} there is a
  subsequence $(v_{n_k})$ such that both $(v_{n_k})$ and $(Sv_{n_k})$
  are equivalent to the unit vector basis of $\ell_1$.
  This contradicts the assumption that $S \in \SS(X)$.
\end{proof}

\section{Applications of Theorem~\ref{main}}

In this section we give applications and corollaries of
Theorem~\ref{main}.

\subsection{}
The first application was obtained by Milman \cite{Milman:70}. By
\cite{KP} we have that for $2< p < \infty$, every weakly null
seminormalized sequence in $L_p[0,1]$ has a subsequence which is
equivalent to the unit vector basis of $\ell_p$ or $\ell_2$.  Thus
$\#\bigl(\SP_{1,w} (L_p[0,1])/\!\!\approx\bigr) =2$. Moreover, $\ell_1
\not \hookrightarrow L_p[0,1]$. Thus, by Theorem~\ref{main}, the
product of any two strictly singular operators on $L_p[0,1]$
($2<p<\infty$) is compact.

\subsection{}
An infinite dimensional subspace $Y$ of a Banach space $X$ is said to
be \term{partially complemented} if there exists an infinite
dimensional subspace $Z\subset X$ such that $Y\cap Z=0$ and $Y+Z$ is
closed. In general, the adjoint of a strictly singular operator
doesn't have to be strictly singular. However, Milman proved
in~\cite{Milman:70} that if $X^*$ is separable and every infinite
dimensional subspace of $X$ is partially complemented, then the
adjoint of every strictly singular operator defined on $X$ is again
strictly singular. Milman then used this fact to show that the product
of any two strictly singular operators on $L_p[0,1]$ ($1<p<2$) is
compact. This can be immediately generalized to the following dual
version of Theorem~\ref{main}.

\begin{corollary}
  Suppose that $X$ is a Banach space such that $X^*$ is separable and
  every infinite dimensional subspace of $X$ is partially
  complemented. If $\bigl(\SP_{1,w}(X^*)/\!\!\approx\bigr)=n$, and
  $T_1,\ldots,T_{n+1}\in\SS(X)$, then $T_{n+1}T_n\dots T_1$ is
  compact. Moreover, if $\ell_1$ does not isomorphically embed in
  $X^*$ then $T_nT_{n-1}\dots T_1$ is compact.
\end{corollary}

\subsection{}
Let $X$ be a reflexive Banach space with a basis $(e_i)$
such that for some $1 \le \xi < \omega_1$ there exists $0<\delta<1$
such that
\begin{equation} \label{asl_1}
  \Bignorm{ \sum_{i=1}^n x_i } \ge
  \delta \sum_{i=1}^n \norm{x_i}
\end{equation}
for any finite block sequence $(x_i)_{i=1}^n$ with $(\min \supp
x_i)_{i=1}^n \in \S_\xi$. Fix $m \in \N$ and by Remark~\ref{R:Sxi}(iv)
let $N= (n_i)$ be a subsequence of $\N$ such that $\S_{\xi m } (N)
\subseteq [ \S_\xi]^m$. Thus for any block sequence $(x_n)$ in $X$ and
any $F \in \S_{\xi m}$ we have
\begin{equation} \label{ell1}
 \Bignorm{ \sum_{i \in F} x_{n_i} } \ge \delta^m \sum_{i\in F}\norm{x_{n_i}}.
\end{equation}
Hence if $(x_{n_i})$ is seminormalized then $(x_{n_i})$ is
$\approx_{\xi m}$-equivalent to the unit vector basis of $\ell_1$.
Therefore the proof of Proposition~\ref{approx}(\ref{app-c}) gives
that if $(x_n)$ is any Schreier spreading sequence in $X$ then $(x_n)$
is $\approx_{\xi m}$-equivalent to the unit vector basis of $\ell_1$.
Since $\ell_1 \not \hookrightarrow X$, by Theorem~\ref{main} we obtain
that $\SS_{\xi m}(X)= \K(X)$. Banach spaces that satisfy
(\ref{asl_1}) are for example Tsirelson type spaces $T[ \delta ,
\S_\xi]$ or more general mixed Tsirelson spaces
$T\bigl[\bigl(\frac{1}{m_i},\S_{n_i}\bigr)_{i \in \N}\bigr]$, or
similar type of hereditarily indecomposable Banach spaces constructed
and studied in~\cite{AD}.

\subsection{}
Let $R$ be the Banach space constructed by C.J.~Read in~\cite{Read}.
It is shown in \cite{Read} that $R$ has precisely two symmetric bases,
(which shall be denoted by $(e_m^Y)_n$ and $(e_n^Z)_n$), up to equivalence.

\begin{proposition}\label{P:R}
  If $(y_n)$ is a Schreier spreading sequence in~$R$, (not necessarily
  symmetric and not necessarily a basis for the whole space), then
  either $(y_n) \approx_1 (e_n^Y)$ or $(y_n) \approx_1 (e_n^Z)$ or
  $(y_n)$ is $\approx_1$-equivalent to the unit vector basis
  of~$\ell_1$.
\end{proposition}

\begin{proof}
  In \cite[page 38, lines 14 and 17]{Read} two norms $\norm{\cdot
  }_Y$ and $\norm{ \cdot }_Z$ are constructed on $c_{00}$ so that the
  standard basis $(e_n)$ of $c_{00}$ is symmetric with respect to
  either norm \cite[page 38, line -2]{Read}. Then $Y$ denotes the
  completion of $(c_{00}, \norm{ \cdot }_Y)$ and $Z$ denotes the
  completion of $(c_{00}, \norm{ \cdot }_Z)$. It is proved in
  \cite[Lemma 2, page 39]{Read} that $Y$ and $Z$ are isomorphic (and
  we denote them by~$R$). Thus if $(e_n^Y)_n$ and $(e_n^Z)_n$ denotes
  the standard basis of $c_{00}$ in $Y$ and $Z$ respectively then
  $(e_n^Y)_n$ and $(e_n^Z)_n$ are normalized symmetric bases for $R$.
  Also, $(e_n^Y)_n$ and $(e_n^Z)_n$ are not equivalent by the
  estimates of \cite[page 39, lines 7 and 9]{Read}. From page~$40$,
  line~$13$ to the end of section~$6$ (page~$47$) it is shown in
  \cite{Read} that if $(y_n)$ is a symmetric $\norm{ \cdot
  }_Y$-normalized block basic sequence of $(e_n^Y)_n$  in $R$
  then $(y_n)$ is equivalent to
  $(e_n^Y)_n$, or $(e_n^Z)_n$, or the unit vector basis of $\ell_1$.
  (Then since $R$ is not isomorphic to $\ell_1$, it is obtained that
  $R$ has exactly two symmetric bases). A closer examination of these
  pages will reveal that it is actually shown that if $(y_n)$ is any
  $\norm{ \cdot }_Y$-normalized block sequence in $R$ then one of the
  following two cases happens:

  \noindent {\em Case 1:} $(y_n)$ has a subsequence $(y_{n_i})_i$
  which is equivalent to the unit vector basis of $\ell_1$ (see
  \cite[page 41, lines 5-7]{Read}). Thus if $(y_n)$ is Schreier
  spreading, then $(y_n)$ is $\approx_1$-equivalent to the unit vector
  basis of $\ell_1$. Moreover, if $(y_n)$ is symmetric (\cite[page 41,
  line 9]{Read}) then $(y_n)$ is equivalent to the unit vector basis
  of $\ell_1$.

  \noindent {\em Case 2:} The limit $\lim_{r \to \infty}
  \bignorm{ \sum_j \lambda_j y_{j+r}}_Y$ is equivalent to either
  $\bignorm{ \sum_j \lambda_j e_j^Y}_Y$ or
  $\bignorm{ \sum_j \lambda_j e_j^Z}_Z$ for every $(\lambda_j) \in
  c_{00}$.

  Indeed, $\bignorm{ \sum_j \lambda_j y_{j+r}}_Y$ is the left hand
  side of the displayed formula \cite[page 47, line 7]{Read} by virtue
  of the notation \cite[page 40, line -4]{Read}.  Thus by
  \cite[page 47, line 7]{Read} the limit $\lim_{r \to \infty}
  \bignorm{ \sum_j \lambda_j y_{j+r}}_Y$ is denoted by $|||
  {\bf \lambda} |||$ in \cite[page 47, line 10]{Read}, or by
  $p({\bf \lambda}, {\bf \beta})$ in \cite[section 7]{Read}.  It is
  concluded in \cite[page 50, line 3]{Read} that $||| {\bf \lambda}
  |||$ is equivalent to either $\bignorm{ \sum_j \lambda_j e_j^Y}_Y$
  or $\bignorm{ \sum_j \lambda_j e_j^Z}_Z$.

  Thus in Case 2, if $(y_n)$ is Schreier spreading then $(y_n)
  \approx_1(e_n^Y)$, or $(y_n) \approx_1(e_n^Z)$. Moreover, if $(y_n)$
  is symmetric \cite[page 47, line 9]{Read} then $(y_n)$ is
  equivalent to $(e_n^Y)$ or $(e_n^Z)$.
\end{proof}

By combining Proposition~\ref{P:R} and Theorem~\ref{main}, we obtain
that the product of any three  operators in $\SS_1(R)$ is
compact.

\subsection{}
Theorem~\ref{main} may also be used to provide invariant subspaces of
operators. A well known theorem of Lomonosov~\cite{Lomonosov:73}
asserts that if $T$ is an operator on a Banach space such that $T$
commutes with a non-zero compact operator, then $T$ has a (proper
non-trivial) invariant subspace.  Moreover, if the Banach space is
over complex scalars and $T$ is not a multiple of the identity, then
there exists a proper non-trivial subspace which is hyperinvariant for
$T$. When the Banach space is over real scalars, one can find a
hyperinvariant subspace for $T$ if $T$ doesn't satisfy an irreducible
quadratic equation, see~\cite{Hooker:81,Sirotkin:05}.

\begin{proposition}\label{invsbsp}
    Suppose that $X$ is a Banach space and $1 \le \xi < \omega_1$.
    \begin{enumerate}
    \item If $\nSPxi$ is finite then every operator $S\in\SS_\xi(X)
      \setminus \{ 0 \}$ has a non-trivial hyperinvariant subspace.
    \item If $\nSP$ is finite then every operator $S\in\SS(X)\setminus
      \{ 0 \}$ has a non-trivial hyperinvariant subspace.
    \end{enumerate}
\end{proposition}

\begin{proof}
  Suppose that either $\nSPxi$ is finite and $S\in\SS_\xi(X)\setminus
  \{ 0 \}$, or $\nSP$ is finite and $S\in\SS(X)\setminus \{ 0 \}$. If
  $S$ has eigenvalues, then every eigenspace is a hyperinvariant
  subspace and we are done.  Suppose $S$ has no eigenvalues. So we can
  assume that $S$ is quasinilpotent with trivial kernel. It follows
  that $S$ doesn't satisfy any real-irreducible quadratic equation.
  Theorem~\ref{main} implies that $S^m$ is compact for some $m$. Also,
  $S^m$ is non-zero as otherwise zero would be an eigenvalue of $S$.
  Since $S$ commutes with $S^m$, it follows that $S$ has a non-trivial
  hyperinvariant subspace.
  \end{proof}

  A similar reasoning shows that if, under the hypotheses of
  Proposition~\ref{invsbsp}, $T$ commutes with $S$ then $T$ commutes
  with the compact operator $S^m$. Therefore, if $S^m\ne0$ and either
  $X$ is a complex Banach space or $X$ is real and $T$ doesn't satisfy
  any irreducible quadratic equation, then $T$ has a hyperinvariant
  subspace.

  Note that Read~\cite{Read2} constructed an example of a strictly
  singular operator with no invariant subspaces.  A further
  application of Proposition~\ref{invsbsp} is
  Corollary~\ref{hyperinv}.

\section{Hereditarily indecomposable Banach spaces}

In \cite{Gowers:93} an infinite dimensional Banach space was defined
to be \term{hereditarily indecomposable (HI)} if for every two
infinite dimensional subspaces $Y$ and $Z$ of $X$ with $Y \cap Z = \{
0 \}$ the projection from $Y + Z$ to $Y$ defined by $y+z \mapsto y$
(for $y \in Y$ and $z \in Z$) is not bounded.
It is is observed in \cite{Gowers:93} that this is equivalent to the
fact that for every two infinite dimensional subspaces $Y$ and $Z$ of
$X$ and for every $\varepsilon >0$ there exists a unit vector $y \in
Y$ such that $\dist \bigl(y, Z \bigr) < \varepsilon$. This motivates
us to introduce the following definition.

\begin{definition}
  Let $1 \le \xi < \omega_1$. We say that a Banach space $X$ is
  $\S_\xi$-hereditary indecomposable (HI$_\xi$) if for every
  $\varepsilon>0$, infinite-dimensional subspace $Y\subseteq X$ and
  basic sequence $(x_n)$ in $X$ there exist an index set $F \in
  \S_\xi$ and a unit vector $y\in Y$ such that the
  $\dist\bigl(y,[x_i]_{i\in F}\bigr)<\varepsilon$.
\end{definition}

It is obvious that if $1 \le \xi < \omega_1$ and $X$ is HI$_{\xi}$
then $X$ is HI. Similarly to Proposition~\ref{P:SS_xi}(\ref{SSb}), if
$X$ is HI$_\xi$ and $\xi<\zeta$ then $X$ is HI$_\zeta$.

\begin{remark}\label{HI_xi}
  Let $X$ be a Banach space and $1 \le \xi < \omega_1$.
  \begin{enumerate}
  \item\label{HIa} $X$ is HI$_\xi$ if and only if for every
    $\varepsilon>0$, infinite-dimensional subspace $Y\subseteq X$ and
    normalized basic sequence $(x_n)$ in $X$ there exist a subsequence
    $(x_{n_k})$, $F \in \S _\xi$ and unit vector $y\in Y$ such that
    $\dist\bigl(y,[x_{n_k}]_{k\in F}\bigr)<\varepsilon$.
  \item\label{HIb} If $X$ is a reflexive Banach space with a basis
    $(e_n)$, then $X$ is HI$_\xi$ if and only if for every
    $\varepsilon>0$, infinite-dimensional block subspace $Y\subseteq
    X$ and normalized block sequence $(y_n)$ of $(e_n)$, there exists
    $G \in \S _\xi$ and unit vector $y \in Y$ such
    that $\dist\bigl(y,[y_i]_{i\in G}\bigr)<\varepsilon$.
  \end{enumerate}
\end{remark}

The proof of (\ref{HIa}) is trivial. For the proof of  (\ref{HIb})
notice that if $X$ is a reflexive Banach space with a basis $(e_n)$, $Y$ is an
infinite dimensional subspace of $X$ and $(x_n)$ is a basic sequence in $X$ then
$\left( \frac{x_n}{\norm{x_n}}\right)$ is weakly null thus by passing to a subsequence and
relabeling we can assume that  $\left( \frac{x_n}{\norm{x_n}}\right)$ is
``near'' a block sequence of $(e_n)$ \cite{BP}. Similarly, $Y$ contains an infinite
dimensional block subspace. The details are left to the reader.

\begin{example}\label{GM}
  \textit{The HI space constructed by Gowers and Maurey \cite{Gowers:93},
  which will be denoted by $GM$, is an HI$_3$ space.}

Indeed, we outline the proof from \cite{Gowers:93} that $GM$ is HI and
we indicate that the proof actually shows that $GM$ is HI$_3$. An
important building block of the proof is the notion of rapidly
increasing sequence vectors (denoted by RIS vectors). Before defining
the RIS vectors, we need to back up and define the $\ell_{1+}^n$
average with constant $1+ \varepsilon$ (for $n \in \N$ and
$\varepsilon >0$). Let $n \in \N$ and $\varepsilon >0$.  We say that a
vector $y \in GM$ is an \term{$\ell_{1+}^n$ average} with constant
$1+\varepsilon$ if $\norm{y}=1$ and $y$ can be written as $y= x_1 +
\cdots +x_n$ where $x_1 < \cdots < x_n$, $x_i$'s are non-zero, and
$\norm{x_i} \le (1+ \varepsilon) n^{-1}$ for every $i$. It is shown in
\cite[Lemma 3]{Gowers:93} that if $U$ is any infinite dimensional
block subspace of $GM$, $\varepsilon >0$ and $n \in \N$ then there
exists $y \in U$ which is an $\ell_{1+}^n$ average with constant $1 +
\varepsilon $. In fact, the proof shows that $(u_n)$ is a block
basis of $U$ then there exists $F \in \S_1$ and $y \in
[u_n]_{n \in F}$ which is an $\ell_{1+}^n$ average with constant $1 +
\varepsilon$.

For $N \in \N$ and $\varepsilon >0$ a vector $z \in GM$
is called an \term{RIS vector} of length $N$ and constant $1+
\varepsilon$ if $z$ can be written as $z=(y_1 + \cdots +y_N)/
\norm{y_1+ \cdots + y_N}$ where $y_1<\dots<y_N$ and each $y_k$ is an
$\ell_{1+}^{n_k}$ average with constant $1+ \varepsilon$ and the
positive integers $(n_k)_{k=1}^N$ are defined inductively to satisfy
$n_1 \ge 4(1+ \varepsilon) 2^{N/\varepsilon'} / \varepsilon'$ where
$\varepsilon ' = \min (\varepsilon , 1)$, and
$\sqrt{ \log_2 (n_{k+1}+1)} \ge 2 \# \supp y_k / \varepsilon'$, where
$\supp y$ stands for the support of the vector $y$ relative to the
standard basis of $GM$. Thus if $N \in \N$,
$\varepsilon >0$ and $U$ is an infinite dimensional block subspace of
$GM$ spanned by a block sequence $(u_n)$, then there exists $G \in
S_2$ and $z\in[u_n]_{N \in G}$ which is an RIS vector of length $N$ and
constant $1+ \varepsilon$.

The idea of the proof that $GM$ is HI is
then the following (\cite[page $868$]{Gowers:93}): Given any $k \in
\N$ and two block subspaces $Y$ and $Z$ of $GM$, spanned by block
sequences $(y_n)$ and $(z_n)$ respectively, let $x_1 \in Y$ be an RIS
of length $M_1:= j_{2k}$ and constant $41/40$ (the sequence $(j_n)$ is
an increasing sequence of integers which is used at the definition of
the space $GM$ \cite[pages $862$ and $863$]{Gowers:93}).  The vector
$x_1$ determines then a positive integer $M_2$. Then a vector $x_2$ is
chosen in $Z$ such that $x_1<x_2$ and $x_2$ is an RIS vector of length
$M_2$ and constant $41/40$.  The vectors $x_1$ and $x_2$ determine a
positive integer $M_3$.  Then a vector $x_2$ is chosen in $Y$ such
that $x_2<x_3$ and $x_3$ is an RIS vector of length $M_3$ and constant
$41/40$.  Continue similarly choosing total of $k$ block vectors $x_i$
alternatingly from $Y$ and $Z$.  Let $y= \sum x_{2i-1}/
\norm{\sum x_{2i-1}} \in Y$ and $z= \sum x_{2i}/\norm{\sum x_{2i}} \in
Z$. Then it is shown that $\norm{y+z} \ge (1/3)\sqrt {\log_2(k+1) }
\norm{y-z}$. Since $k$ is arbitrary, this shows that $GM$ is HI. By
the remarks about the support of an RIS vector, one can make sure that
there exist $H_1, H_2 \in \S_3$ such that $y \in [y_n]_{n \in H_1}$
and $z \in [z_n]_{n \in H_2}$. This proves that $GM$ is an HI$_3$
space. Proposition~\ref{OperatorsonHI} implies that if $GM$ is considered
as a complex Banach space then every operator on $GM$ can be written
in the form $\lambda + S$ where $\lambda \in \C$ and $S \in \SS _3(GM)$.
\end{example}

\begin{example}\label{AD}
  \textit{The HI space constructed by S.A.~Argyros and I.~Deliyanni \cite{AD},
  which will be denoted by $AD$, is an HI$_{\omega 3}$ space.}

  This can be done similarly to the Example~\ref{GM} by closely
  examining the proof of \cite{AD} showing that $AD$ is an HI space.
\end{example}

Next we use Desriptive Set Theory in order to prove the following
result which signifies the importance of separable $\S_\xi$-HI Banach
spaces and $\S_\xi$-strictly singular operators defined on separable
Banach spaces.

\begin{theorem}\label{t1}
  Let $X, Y$ be separable Banach spaces and $S\in\mathcal{L}(X, Y)$.
  Then the following hold.
  \begin{enumerate}
    \item\label{a} $X$ is HI if and only if $X$ is HI$_{\xi}$ for some $\xi<\omega_1$.
    \item\label{b} $S$ is strictly singular if and only if $S$ is
      $\S_\xi$-strictly singular for some $\xi<\omega_1$.
  \end{enumerate}
\end{theorem}

For the proof of Theorem \ref{t1} we need some results from
Descriptive Set Theory which we briefly recall.  \bigskip

\noindent \textit{Trees.} Let $\bt$ by the set of all finite
sequences of natural numbers. By $\ibbt$ we shall denote the
subset of $\bt$ consisting of all strictly increasing finite
sequences. We view $\bt$ as a tree equipped with the (strict)
partial order $\sqsubset$ of extension. A \textit{tree} $T$ on
$\N$ is a downwards closed subset of $\bt$. By $\tr$ we denote
the set of all trees on $\N$. Thus
\begin{displaymath}
  T\in\tr \Leftrightarrow \forall s,t\in\bt \
  (s\sqsubset t \text{ and } t\in T \Rightarrow s\in T).
\end{displaymath}
Notice that $\ibbt$ belongs to $\tr$.

By identifying every $T\in\tr$ with its characteristic function (i.e.
an element of $2^{\bt}$), it is easy to see that the set $\tr$ becomes
a closed subset of $2^{\bt}$.  For every $\sigma\in \N^\N$ and every
$n\in\N$ we let $\sigma|n=\bigl(\sigma(1),\dots,\sigma(n)\bigr)\in
\bt$.  A tree $T\in\tr$ is said to be \textit{well-founded} if for
every $\sigma\in \N^\N$ there exists $n\in\N$ such that
$\sigma|n\notin T$. By $\wf$ we denote the subset of $\tr$ consisting
of all well-founded trees.

For every $T\in\tr$ we let
$T'=\bigl\{s\in T\mid \exists t\in T \text{ with } s\sqsubset t\bigr\}$.
Observe that $T'\in\tr$. By transfinite recursion, for every $T\in\tr$
we define $(T^{(\xi)})_{\xi<\omega_1}$ as follows. We set $T^{(0)}=T$,
$T^{(\xi+1)}= \big( T^{(\xi)}\big)'$ and
$T^{(\lambda)}=\bigcap_{\xi<\lambda} T^{(\xi)}$ if $\lambda$ is limit.
Notice that $T\in\wf$ if and only if the sequence
$(T^{(\xi)})_{\xi<\omega_1}$ is eventually empty. For every $T\in\wf$,
the \textit{order} $o(T)$ of $T$ is defined to be the least countable
ordinal $\xi$ such that $T^{(\xi)}=\varnothing$. We will need the
following Boundedness Principle for $\wf$, \cite[Theorem~31.2]{Kechris}.

\begin{theorem}
\label{bwf} If $A$ is an analytic subset of $\wf$, then
$\sup\{ o(T)\mid T\in A\}<\omega_1$.
\end{theorem}

If $S, T\in \tr$, then a map $\phi:S\to T$ is said to
be \textit{monotone} if for every $s_1, s_2\in S$ with
$s_1\sqsubset s_2$ we have $\phi(s_1)\sqsubset \phi(s_2)$.
Notice that if $S, T$ are well-founded trees and
there exists a monotone map $\phi:S\to T$, then
$o(S)\le o(T)$. Also observe that for every $\xi<\omega_1$
the Schreier family $\S_\xi$ is a well-founded tree
and $o(\S_\xi)\ge \xi$.
\bigskip

\noindent \textit{Standard Borel spaces.} Let $(X,\Sigma)$ be a
measurable space, i.e. $X$ is a set and $\Sigma$ is a $\sigma$-algebra
on $X$. The pair $(X,\Sigma)$ is said to be a
\textit{standard Borel space} if there exists a Polish topology $\tau$
on $X$ such that the Borel $\sigma$-algebra of $(X,\tau)$ coincides
with~$\Sigma$.  Invoking the classical fact that for every Borel
subset $B$ of a Polish space $(X,\tau)$ there exists a finer Polish
topology $\tau'$ on $X$ making $B$ clopen and having the same Borel
set as $(X,\tau)$ (see \cite[Theorem 13.1]{Kechris}), we see that if
$(X,\Sigma)$ is a standard Borel space and $B\in \Sigma$, then $B$
equipped with the relative $\sigma$-algebra is a standard Borel space
too.

Let $X$ be a Polish space and denote by $F(X)$ the set of all closed
subsets of~$X$. We endow $F(X)$ with the $\sigma$-algebra $\Sigma$
generated by the sets
$$\bigl\{ F\in F(X)\mid F\cap U\neq\varnothing\bigr\},$$
where $U$ ranges over all non-empty open subsets of
$X$. The measurable space $\bigl(F(X),\Sigma\bigr)$ is called
the Effros-Borel space of~$X$. It is well-known
that the Effros Borel space is a standard Borel space,
see \cite[ Theorem 12.6]{Kechris}.

Now let $X$ be a separable Banach space. Denote by $\subs(X)$ the set
of all infinite-dimensional subspaces of $X$. It is easy to see that
$\subs(X)$ is a Borel subset of $F(X)$ (see \cite{Kechris}, Exercises
12.19 and 12.20), and so, a standard Borel space on its own right.  We
will need the following fact, which was isolated explicitly in
\cite{ADo}. Its proof follows by a straightforward application of the
Kuratowski--Ryll-Nardzewski selection theorem see
\cite[Theorem 12.13]{Kechris}.

\begin{proposition}\label{sel}
  Let $X$ be a separable Banach space.  There exists a sequence
  $S_l:\subs(X)\to X$, $l\in\N$, of Borel maps such that for every
  subspace $Y$ of $X$ the sequence $\big(S_l(Y)\big)$ is in the sphere
  $S_Y$ of $Y$ and, moreover, it is norm dense in~$S_Y$.
\end{proposition}

For more background material on $\subs(X)$ we refer
to \cite{ADo}, \cite{B} and \cite{Kechris}. We are ready to
proceed to the proof of Theorem \ref{t1}.

\begin{proof}[Proof of Theorem \ref{t1}]
  \eqref{a} Clearly we only need to show that if $X$ is HI, then $X$
  is HI$_{\xi}$ for some $\xi<\omega_1$.  So, fix a separable HI
  Banach space~$X$. Let
  $$\B=\bigl\{ (x_n)\in X^\N\mid (x_n) \text{ is a normalized
    basic sequence in }
  X\bigr\}.$$
  We claim that $\B$ is $F_\sigma$ in $X^\N$, the later equipped with
  the product topology. To see this, for every $k\in\N$ let $\B_k$ be
  the set of all normalized basic sequences $(x_n)$ with basis
  constant less or equal than $k$. It is easy to see that $\B_k$ is
  closed in $X^\N$.  As $\B$ is the union over all $k\in\N$ of $\B_k$,
  this shows that $\B$ is~$F_{\sigma}$.  Since $X$ is separable,
  $X^\N$ is Polish. Thus $\B$ is a standard Borel space.

  Let $S_l:\subs(X)\to X$, $l\in\N$, be the sequence of Borel maps
  obtained by Proposition~\ref{sel}. For every $m\in\N$, every
  $(x_n)\in \B$ and every $Y\in\subs(X)$ we define a tree $T=\Tm\in
  \tr$ to be the set of all $t=(l_1<\dots< l_k)\in\ibbt$ such that
  \begin{displaymath}
    \Bignorm{S_l(Y)- \sum_{i=1}^k a_i x_{l_i}}\ge\tfrac{1}{m}\quad
    \text{for any }a_1,\dots,a_k\mbox{ in }\mathbb Q
    \text{ and any $l$ in $\N$}.
  \end{displaymath}
  For every $m\in\N$ consider the map $\Phi_m:\B\times \subs(X)\to
  \tr$ defined by
  $$\Phi_m\bigl((x_n),Y\bigr)=\Tm.$$
  \noindent \textsc{Claim 1.} \textit{The following hold.
  \begin{enumerate}
   \item\label{Borel} For every $m\in\N$ the map $\Phi_m$ is Borel.
   \item\label{wf} For every $m\in\N$, every $(x_n)\in\B$ and every
    $Y\in\subs(X)$ the tree $\Tm$ is well-founded.
   \item\label{est} Let $\zeta<\omega_1$ and assume that $X$ is not
    HI$_{\zeta}$. Then there exist $m\in\N$, $(x_n)\in\B$
    and $Y\in\subs(X)$ such that $o\Bigl(\Tm\Bigr)\ge\zeta$.
  \end{enumerate} }
  \bigskip

  \noindent \textit{Proof of the claim.} \eqref{Borel} For those
  readers familiar with descriptive set theoretic computations, this
  part of the claim is a straightforward consequence of the definition
  of the tree $\Tm$. However, for the convenience of the readers not
  familiar with these computations, we shall describe a more detailed
  argument.

  Fix $m\in\N$. For every $t\in\bt$ let $U_t=\{ T\in\tr\mid t\in T\}$.
  As the topology on $\tr$ is the pointwise one, we see that the
  family $\{U_t\mid t\in \bt\}$ forms a sub-basis of the topology on
  $\tr$. Thus, it is enough to show that for every $t\in\bt$ the set
  $$\Phi_m^{-1}(U_t)=\Bigl\{\bigl((x_n),Y\bigr)\in \B\times
  \subs(X)\mid t\in\Tm\Bigr\}$$
  is Borel. So, let $t\in \bt$.
  If $t\notin \ibbt$, then $\Phi_m^{-1}(U_t)=\varnothing$. Hence, we
  may assume that $t=(l_1<\dots<l_k)\in \ibbt$.

  For every $j\in\N$ the map $\pi_j: \B\times \subs(X)\to X$ defined
  by $\pi_j\bigl((x_n), Y\bigr)=x_j$ is clearly Borel.  For every
  $\mathbf{a}=(a_i)_{i=1}^k\in\mathbb{Q}^k$ and every $l\in\N$
  consider the map $H_{\mathbf{a},l}:\B\times \subs(X)\to\mathbb R$
  defined by
  \begin{displaymath}
    H_{\mathbf{a},l}\bigl((x_n),Y\bigr)=
    \Bignorm{S_l(Y)-\sum_{i=1}^k a_i x_{l_i}}=
    \Bignorm{S_l(Y)-\sum_{i=1}^k a_i \pi_{l_i}\bigl((x_n),Y\bigr)}.
  \end{displaymath}
Invoking the above remarks, the Borelness of the map $S_l$ and the
continuity of the norm, we see that the map $H_{\mathbf{a},l}$ is Borel.
Thus, setting
\begin{math}
  A_{\mathbf{a},l,m}= H^{-1}_{\mathbf{a},l}\bigl(
  [\tfrac{1}{m},+\infty)\bigr)
\end{math}
we get that $A_{\mathbf{a},l,m}$ is a Borel subset of $\B\times \subs(X)$
for every $\mathbf{a}\in\mathbb{Q}^k$, every $l\in\N$, and every $m\in\N$.
It follows from
\begin{displaymath}
  \Phi_m^{-1}(U_t)= \bigcap_{\mathbf{a}\in\mathbb{Q}^k}
  \bigcap_{l\in\N} A_{\mathbf{a},l,m}
\end{displaymath}
that $\Phi_m^{-1}(U_t)$ is a Borel subset
of $\B\times \subs(X)$, as desired.

\noindent\eqref{wf} Assume,  towards a contradiction, that there exist
$m\in\N$, $(x_n)\in\B$ and $Y\in\subs(X)$ such that the tree $\Tm$ is
not well-founded. Thus, there exists $\sigma\in\N^\N$ such that
$\sigma|k\in\Tm$ for all $k\in\N$.  Set $n_k=\sigma(k)$. Notice that
$n_k<n_{k+1}$ for every $k\in\N$. Let $Z=[x_{n_k}]$. As the sequence
$\bigl(S_l(Y)\bigr)$ is norm dense in $S_Y$, we see that
$\mathrm{dist}(y,Z)\ge \frac{1}{m}$ for every
$y\in S_Y$. Thus $X$ is not HI, a contradiction.

\noindent\eqref{est} Let $\zeta<\omega_1$ such that $X$ is not HI$_{\zeta}$.  By
definition, there exist $\varepsilon>0$, $(x_n)\in\B$ and $Y\in
\subs(X)$ such that for every $y\in S_Y$ and every $F\in\S_{\zeta}$ we
have $\mathrm{dist}(y,[x_n]_{n\in F})\ge\varepsilon$. Let $m\in\N$
with $\frac{1}{m}<\varepsilon$. It follows that for every
$F=\{l_1<\dots< l_k\}\in\S_\zeta$ we have $F\in\Tm$. Hence, the
identity map $\mathrm{Id}:\S_\zeta\to\Tm$ is a well-defined monotone
map, and so $o\Bigl(\Tm\Bigr)\ge o(\S_\zeta)\ge \zeta$.  The claim is
proved.

\bigskip

We set
\begin{displaymath}
  A=\bigcup_{m\in\N} \Phi_m\bigl(\B\times\subs(X)\bigr)=
  \Bigl\{\Tm\mid m\in\N, (x_n)\in\B \text{ and }
  Y\in\subs(X)\Bigr\}.
\end{displaymath}
By Claim~1\eqref{Borel}, we see that $A$ is an analytic
subset of $\tr$. By Claim~1\eqref{wf}, we get that $A\subseteq \wf$.
Hence, by Theorem~\ref{bwf}, there exists a countable ordinal
$\xi$ such that
\[ \sup\{ o(T):T\in A\}<\xi.\]
Finally, by Claim~1\eqref{est}, we conclude that $X$ is
HI$_{\xi}$, as desired.

\medskip

\eqref{b} The proof is similar to that of part (a). Again it is
enough to show that if $X, Y$ are separable Banach space
and $S\in\L(X,Y)$ is strictly singular, then $S$
is $\S_\xi$-strictly singular for some $\xi<\omega_1$.
As in the previous part, let $\B$ be the $F_\sigma$ subset
of $X^\N$ consisting of all normalized basic sequences in $X$.
For every $m\in\N$ and every $(x_n)\in\B$ we define a tree
$T=T\bigl(m,(x_n)\bigr)\in\tr$ to be the set of all
$t=(l_1< ... < l_k)\in\ibbt$ such that
\begin{displaymath}
  \Bignorm{S\bigl(\sum_{i=1}^k a_i x_{l_i}\bigr)} \ge
  \tfrac{1}{m}\Bignorm{\sum_{i=1}^k a_i x_{l_i}}
\end{displaymath}
for any $a_1,\dots,a_k$ on $\mathbb Q$.
For every $m\in\N$ consider the map $\Psi_m:\B\to \tr$
defined by
$$\Psi_m\big( (x_n)\big)= T\bigl(m,(x_n)\bigr).$$
We have the following analogue of Claim~1. The proof
is identical and it is left to the reader.
\medskip

\noindent \textsc{Claim 2.} The following hold.
\begin{enumerate}
\item\label{2-Borel} For every $m\in\N$ the map $\Psi_m$ is Borel.
\item\label{2-wf} For every $m\in\N$ and every $(x_n)\in\B$ the tree
  $T\bigl(m,(x_n)\bigr)$ is well-founded.
\item\label{2-est} Let $\zeta<\omega_1$ and assume that $S$ is not
  $\S_{\zeta}$-strictly singular. Then there exist $m\in\N$
  and $(x_n)\in\B$ such that $o\bigl(T(m,(x_n))\bigr)\ge \zeta$.
\end{enumerate}
\bigskip

By Claim~2\eqref{2-Borel} and \eqref{2-wf}, and Theorem \ref{bwf},
there exists a countable ordinal $\xi$ such that
$$\sup\Bigl\{o\bigl(T(m,(x_n))\bigr)\mid m\in\N \text{ and }
  (x_n)\in\Bigr\}< \xi.$$
Hence, by Claim 2\eqref{2-est}, we conclude that $S$ is
$\S_\xi$-strictly singular. The proof is completed.
\end{proof}

\begin{remark}
  We notice that part~\eqref{b} of Theorem~\ref{t1} is not valid if
  $X$ and $Y$ are non-separable. To see this, for every $\xi<\omega_1$
  let $X_\xi$ and $Y_\xi$ be separable Banach spaces and $T_\xi\in
  \L(X_\xi, Y_\xi)$ be a strictly singular operator which is
  not $\S_\xi$-strictly singular and with $\norm{T_\xi}=1$. We let
  $$X=\Big( \sum_{\xi<\omega_1} \oplus X_\xi\Big)_{\ell_1}
  \ \text{ and } \ Y=\Big( \sum_{\xi<\omega_1}\oplus
  Y_\xi\Big)_{\ell_2}.$$
  One can easily ``glue" the sequence $(T_\xi)_{\xi<\omega_1}$ to
  produce a strictly singular operator $T\in\L(X,Y)$ which is not
  $\mathcal{S}_\xi$-strictly singular for any $\xi<\omega_1$.
\end{remark}

In their celebrated paper \cite{Gowers:93}, W.T.~Gowers and B.~Maurey
showed that if $X$ is a complex HI Banach space then every operator $T
\in \L (X)$ can be written as a strictly singular perturbation of a
scalar operator. The proof is based on the definition of the infinite
singular values of an operator and an important fact that is proved
about them.  We recall the definition: Let $X$ be a complex Banach
space, $T \in \L (X)$. We
say that $T$ is \term{infinitely singular} if no restriction of $T$ to
a subspace of finite codimension is an isomorphism.

\begin{lemma}[\cite{Gowers:93}] \label{infsing} If $X$ is an infinite
  dimensional Banach space over $\C$ and $T \in \L (X)$ then there
  exists
  $\lambda \in \C$ such that $T-\lambda I$ is infinitely singular.
\end{lemma}

Using this fact Gowers and Maurey proved the following.

\begin{theorem}[\cite{Gowers:93}] \label{GM-SS}
  Every operator on a complex HI space is of the form $\lambda I+S$
  where $\lambda \in \C$ and $S$ is strictly singular.
\end{theorem}

We use Lemma~\ref{infsing} in the proof of the following result which is
analogous to Theorem~\ref{GM-SS}.

\begin{proposition} \label{OperatorsonHI} If $1 \le \xi < \omega_1$
  and $X$ is a complex HI$_\xi$ space then every $T \in \L (X)$ can be
  written as $T = \lambda I + S$ where $\lambda \in \C$ and $S \in
  \SS_\xi (X)$.
\end{proposition}

\begin{proof}
  Let $X$ be a complex HI$_\xi$ space and $T \in \L (X)$. Assume that
  $T$ is not a scalar multiple of the identity, else there is nothing
  to prove.  By Lemma~\ref{infsing} there exists $\lambda \in \C$ such
  that $S=T-\lambda I$ is infinitely singular. We will show that
  $S\in\SS_\xi (X)$. Let $(x_n)$ be a
  normalized basic sequence in $X$ and $\varepsilon>0$.
  Proposition~2.c.4 of~\cite{LT1} asserts that there is an infinite
  dimensional subspace $Y$ of $X$ such that
  $\norm{S_{|Y}}<\frac{\varepsilon}{3}$. Since $X$ is HI$_\xi$ there
  exists $F \in \S_\xi$, a unit vector $y \in Y$ and a vector
  $x\in[x_n]_{n\in F}$ such that $\norm{y-x}<\frac{\varepsilon}{3 \norm{S} +\varepsilon}$.
  It can then be easily checked that
  $\Bignorm{ \frac{x}{\norm{x}} - x } < \frac{
    \varepsilon}{3 \norm{S}}$. Hence
  \begin{multline*}
    \Bignorm{ S \frac{x}{\norm{x}}} \le \norm{S y } + \norm{S(y-x)} +
    \Bignorm{S \Bigl( x - \frac{x}{\norm{x}} \Bigr) } \\
    \le \frac{\varepsilon}{3} + \norm{S} \frac{\varepsilon}{3 \norm{S}} +
    \norm{S} \Bignorm{x - \frac{x}{\norm{x}}} < \varepsilon.
  \end{multline*}
  Since $\frac{x}{\norm{x}} \in [x_n]_{ n \in F}$, we obtain that $S
  \in \SS_\xi (X)$, which finishes the proof.
\end{proof}

Propositions~\ref{OperatorsonHI} and~\ref{invsbsp} yield the
following result.

\begin{corollary} \label{hyperinv} If $X$ is an infinite dimensional
  complex HI$_\xi$ Banach space for some $1 \le \xi < \omega_1$, such
  that $\nSPxi < \infty$, then every operator $T
  \in \L (X)$ which is not a multiple of the identity has a
  non-trivial hyperinvariant subspace.
\end{corollary}

\begin{question}
  Does there exist any Banach space which satisfies the assumptions of
  Corollary~\ref{hyperinv}?
\end{question}

Finally we examine operators originating from a subspace of an HI$_\xi$
Banach space $X$ and taking values in $X$. The next result will be important in their study:

 \begin{theorem}\label{no-xi}
  If $1\le\xi<\omega_1$ and $X$ is a HI$_\xi$ Banach space then
$\SS_\xi(X,Y)=\SS(X,Y)$ for every Banach space $Y$.
\end{theorem}

\begin{proof}
  It follows from Proposition~\ref{P:SS_xi}(\ref{SSa}) that
  $\SS_\xi(X,Y)\subseteq\SS(X,Y)$. Let $T\in\SS(X,Y)$, $(x_n)$ a basic
  sequence in $X$, and $0<\varepsilon<1$. Choose $\delta>0$ such that
  $\frac{\delta(1+\norm{T})}{1-\delta}<\varepsilon$.  By
  Proposition~2.c.4 of~\cite{LT1} there is an infinite dimensional
  subspace $Z\subseteq X$ such that $\norm{T_{|Z}}<\delta$. Since $X$
  is HI$_\xi$, there exists $F\in\S_\xi$ and vectors $x\in[x_n]_{n\in
    F}$ and $z\in Z$ such that $\norm{z}=1$ and
  $\norm{x-z}<\delta$. It follows that $\norm{x}>1-\delta$ and
  $$\norm{Tx}\le\norm{Tz}+\norm{T}\norm{x-z}<
    \delta\bigl(1+\norm{T}\bigr)<\varepsilon\norm{x}.$$
\end{proof}

We now extend the following result of V.~Ferenczi  (which in turn
is a generalization of Theorem~\ref{GM-SS}).

\begin{theorem}[\cite{F}] \label{Fer}
  If $X$ is a complex HI Banach space, $Y$ is an infinite dimensional
  subspace of $X$ and $T \in \L (Y, X)$ then there exists
  $\lambda\in\C$ and $S \in \SS(Y,X)$
  such that $T= \lambda i_{Y,X} + S$ where $i_{Y,X}: Y \to X$ is the
  inclusion map.
\end{theorem}

\begin{corollary} \label{T:F}  If $X$ is a complex HI$_\xi$ Banach
  space for some $1 \le \xi < \omega_1$, $Y$ is an infinite
  dimensional subspace of $X$ and $T \in \L (Y,X)$, then there exists
  $\lambda\in\C$ and $S \in \SS_\xi (Y,X)$ such that $T= \lambda
  i_{Y,X} + S$.
\end{corollary}

\begin{proof}
  It follows from Theorem~\ref{Fer} that $T=\lambda i_{Y,X}+S$ for
  some $\lambda\in\C$ and $S\in\SS(Y,X)$. Since a subspace of an
  HI$_\xi$-space is again an HI$_\xi$-space, Theorem~\ref{no-xi}
  yields that $S\in\SS_\xi(Y,X)$.
\end{proof}

\end{document}